\documentclass[12pt]{article}
\usepackage{amssymb}
\usepackage{amsfonts}
\usepackage{latexsym}
\usepackage{amsmath}
\usepackage{epsfig}
\usepackage{graphicx,float}

\author{Paul Blanchard \\
F\.{i}gen \c{C}\.{i}l\.{i}ng\.{i}r\thanks{The second 
author would like to thank the
Department of 
Mathematics and Statistics at Boston University for its hospitality 
while this work was in progress.
In addition, she would also like to thank TUB\.{I}TAK for their
support while this research was in progress.} \\
Daniel Cuzzocreo\\
Robert L.~Devaney \\
Daniel M.~Look \\
Elizabeth D.~Russell
}

\title{Checkerboard Julia Sets for Rational Maps} 

\date{Current version: June 16, 2011\\
Original draft: February 8, 2010}

\begin{document}
\catcode`@=11 \font\tenmsy=msbm10 scaled \magstephalf
\font\sevenmsy=msbm8 \font\fivemsy=msbm6
\newfam\msyfam
\textfont\msyfam=\tenmsy \scriptfont\msyfam=\sevenmsy
\scriptscriptfont\msyfam=\fivemsy
\def\msy@{\hexnumber@\msyfam}
\def\Bbb{\ifmmode\let\next\Bbb@\else
\def\next{\errmessage{Use \string\Bbb\space only in math
mode}}\fi\next}
\def\Bbb@#1{{\Bbb@@{#1}}}
\def\Bbb@@#1{\fam\msyfam#1}
\catcode`@=12
\newcommand{\bbZ}{{\Bbb Z}}
\newcommand{\bbQ}{{\Bbb Q}}
\newcommand{\bbR}{{\Bbb R}}
\newcommand{\bbD}{{\Bbb D}}
\newcommand{\bbC}{{\Bbb C}}
\newcommand{\bbN}{{\Bbb N}}
\newcommand{\bbS}{{\Bbb S}}
\newcommand{\be}{\partial B_\lambda}
\newcommand{\B}{B_\lambda}
\newcommand{\T}{T_\lambda}
\newcommand{\Hla}{H_\lambda}
\newcommand{\bla}{\beta_\lambda}
\newcommand{\tla}{\tau_\lambda}
\newcommand{\F}{F_\lambda}
\newcommand{\G}{G_\lambda}
\newcommand{\eps}{\epsilon}
\newcommand{\la}{\lambda}
\newcommand{\Arg}{{\rm Arg}\,}
\newcommand{\ov}{\overline}
\newcommand{\id}{\thicksim}
\baselineskip = 17pt
\newtheorem{definition}{Definition}[section]
\newtheorem{theorem}[definition]{Theorem}
\newtheorem{corollary}[definition]{Corollary}
\newtheorem{proposition}[definition]{Proposition}
\newtheorem{lemma}[definition]{Lemma}
\begin{titlepage}
\maketitle
\noindent {\bf{Abstract.}} In this paper, we consider the family of
rational maps 
$$\F(z) = z^n + \frac{\la}{z^d},$$ 
where $n \geq 2$, $d\geq 1$, and
$\la \in \bbC$.  We consider the case where
 $\la$ lies in the main cardioid of one of the $n-1$ 
principal Mandelbrot sets in these families.  We show that the
Julia sets of these maps are always homeomorphic.
However, two such
maps $\F$ and $F_\mu$ are
conjugate on these Julia sets only if the parameters
at the centers of the given cardioids
satisfy $\mu = \nu^{ j(d+1)}\la$ or 
$\mu = \nu^{j(d+1)}\overline{\la}$
where 
$j \in \bbZ$ and  $\nu$ is an
$n-1^{\rm st}$ root of unity.  We define
a dynamical invariant, which we call  the minimal rotation number. It determines which of these 
maps are conjugate on their Julia sets, and we obtain an exact count of the number of distinct conjugacy classes of maps
drawn from these main cardioids. 

\end{titlepage} 

In recent years there have been many papers dealing with the family of
rational maps given by
$$
\F(z) = z^n + \frac{\la}{z^d},
$$
where $n \geq 2$, $d\geq 1$, and $\la \in \bbC$ \cite{D5}.   For many parameter
values, the Julia sets for these maps are Sierpi\'nski
curves, i.e., planar sets that are homeomorphic to
the well-known Sierpi\'nski carpet fractal. One distinguishing
property of Sierpi\'nski curve Julia sets is that the Fatou set
consists of infinitely many
open disks, each bounded by a simple closed curve,
but no two of these bounding curves intersect.
 
There are many
different ways in which these Sierpi\'nski curves arise as Julia sets
in these families.
For example, the Julia set is a Sierpi\'nski curve 
if $\la$ is a parameter for which
\begin{enumerate}
\item the critical orbits  enter the immediate basin of attraction of
$\infty$ after two or more iterations \cite{DLU};
\item the parameter lies in the main cardioid of a ``buried'' baby
Mandelbrot set \cite{DL1}; or
\item the parameter lies on a buried point in a Cantor necklace in the
parameter plane \cite{D3}.
\end{enumerate}

The parameter planes for these
maps in the cases where $n=d=3$ and $n=d=4$ are shown in Figure~1.  
The red disks 
not centered at the origin are regions where the first case above
occurs. These disks are called Sierpi\'nski holes. 

Many Mandelbrot sets are visible in Figure~1. 
The ones that touch the external red region are not ``buried,'' so their
main cardioids do not contain Sierpi\'nski curve Julia sets. Only the
ones that do not meet this boundary contain parameters from case~2.

Finally, numerous Cantor necklaces, i.e., sets
homeomorphic to the Cantor middle-thirds set with the removed open
intervals replaced by open disks, appear in these figures. The buried
points in the Cantor set portion of the necklace are the parameters
for which case~3 occurs.

\begin{figure}[ht]
\begin{center}
{\includegraphics[height=2.5in]{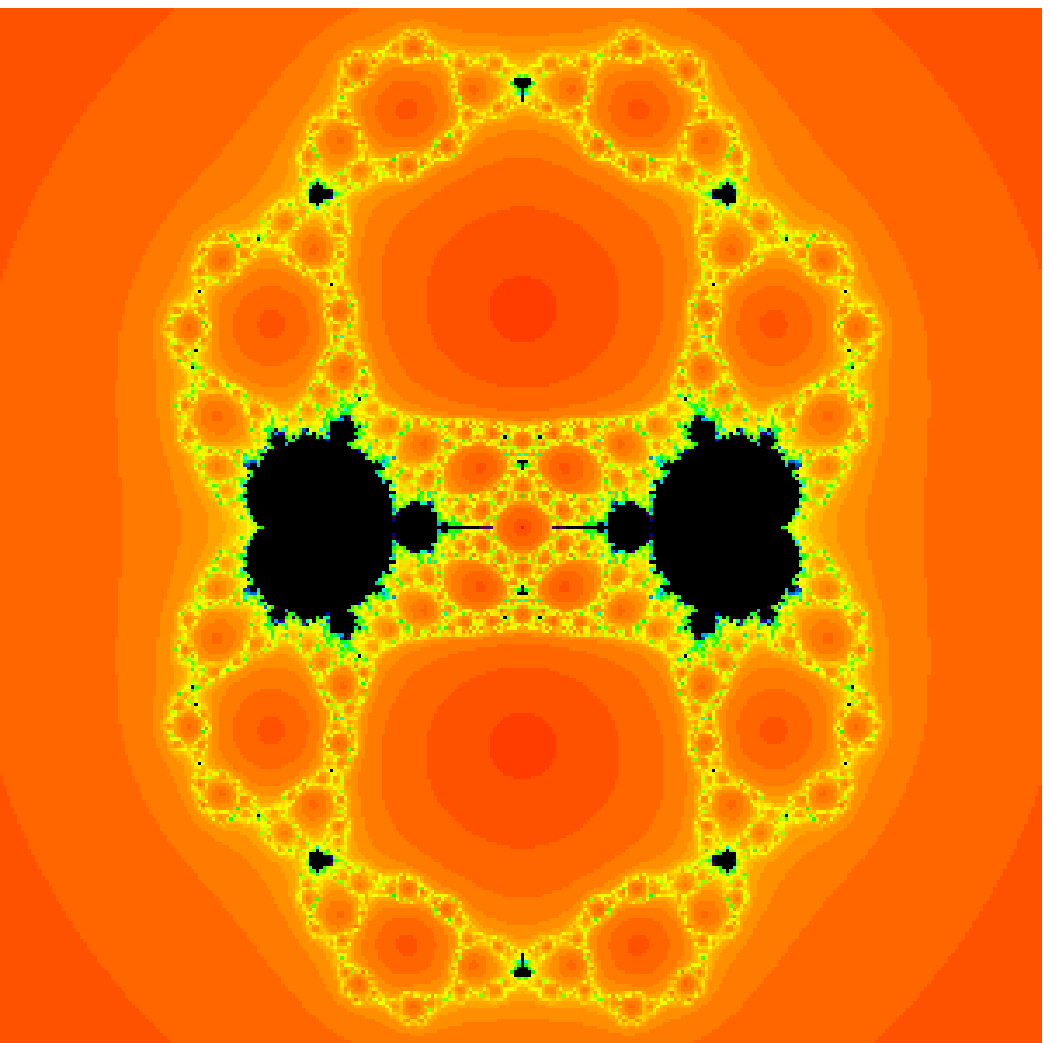}\quad 
\includegraphics[height=2.5in]{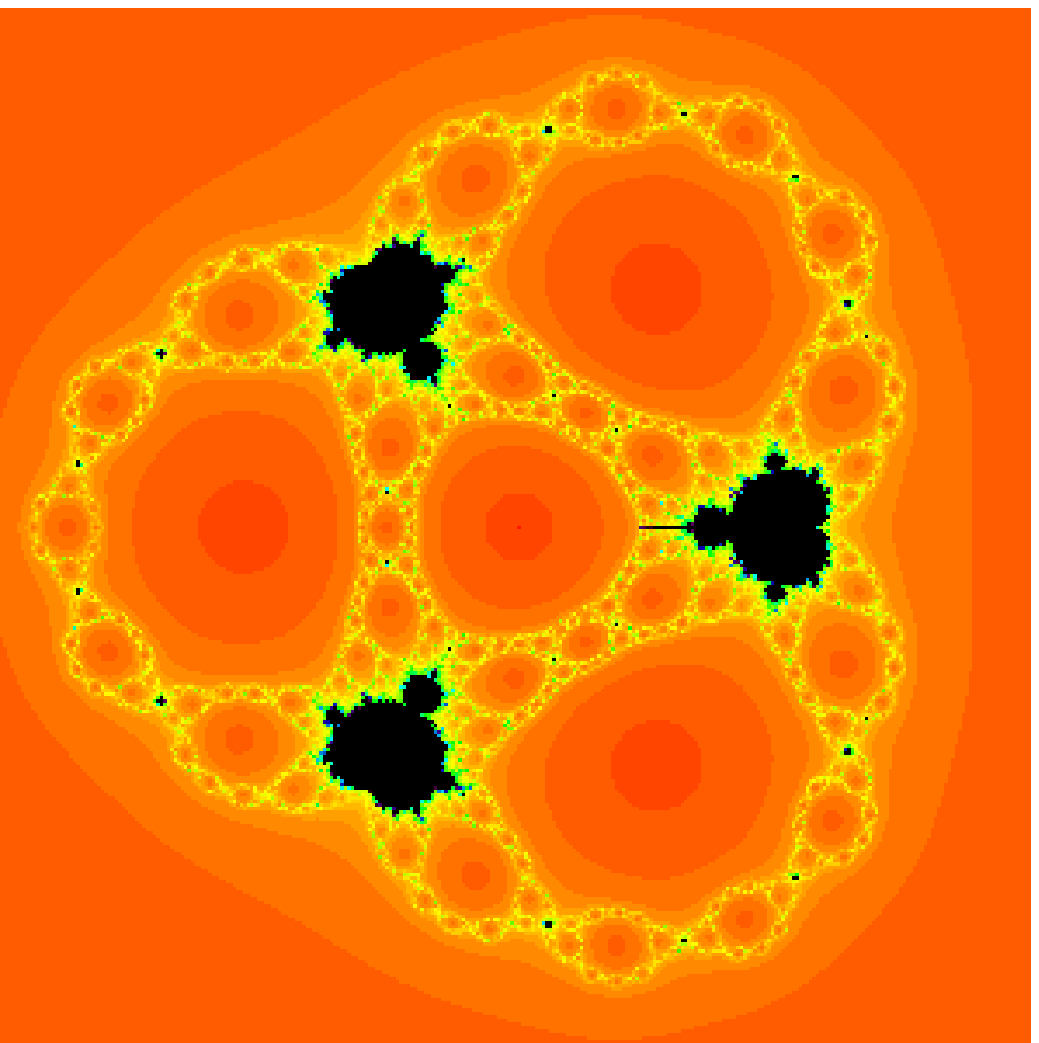}}
\caption{Two parameter planes:  $n=d=3$ (left) and $n=d=4$ (right).}
\label{pplanes}
\end{center}
\end{figure}
 
The dynamical behavior on Sierpi\'nski
curve Julia sets drawn from non-symmetrically located 
Sierpi\'nski holes is never the same \cite{DP}.
That is, only symmetrically located Sierpi\'nski holes
contain parameters for which the corresponding maps have
conjugate dynamics.  While it is
known that two such maps are not conjugate on their Julia sets, there
is no known dynamical invariant that explains this lack of conjugacy. 

In this paper, we describe the topology of and dynamics on
a very different type
of Julia set that arises in these families, the ``checkerboard'' Julia
sets. The parameter spaces contain $n-1$ ``principle Mandelbrot sets," the largest 
Mandelbrot sets in the parameter space~\cite{D4}. We consider the Julia sets for the   
parameters that lie in the main
cardioids of these Mandelbrot sets.

\begin{figure}[ht]
\begin{center}
{\includegraphics[height=2.5in]{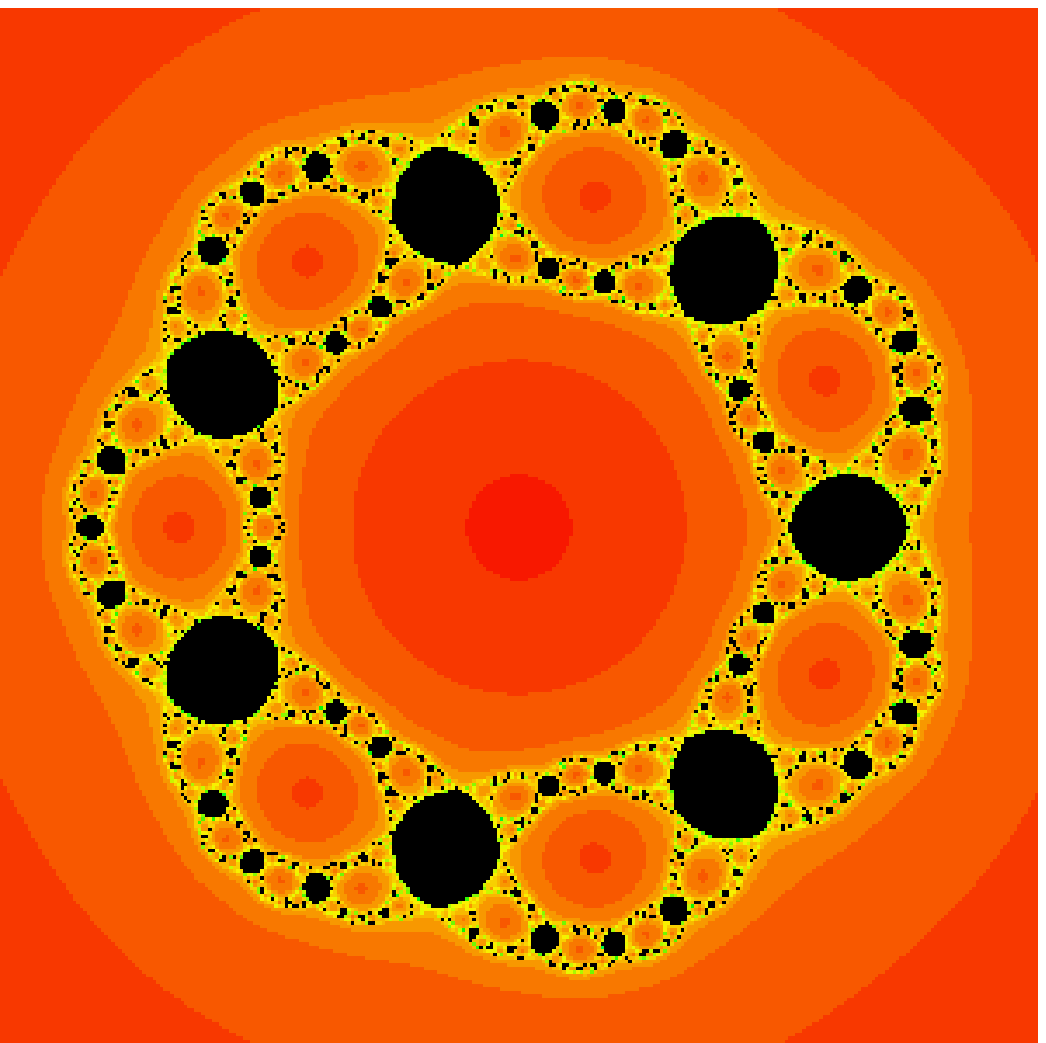}\quad 
\includegraphics[height=2.5in]{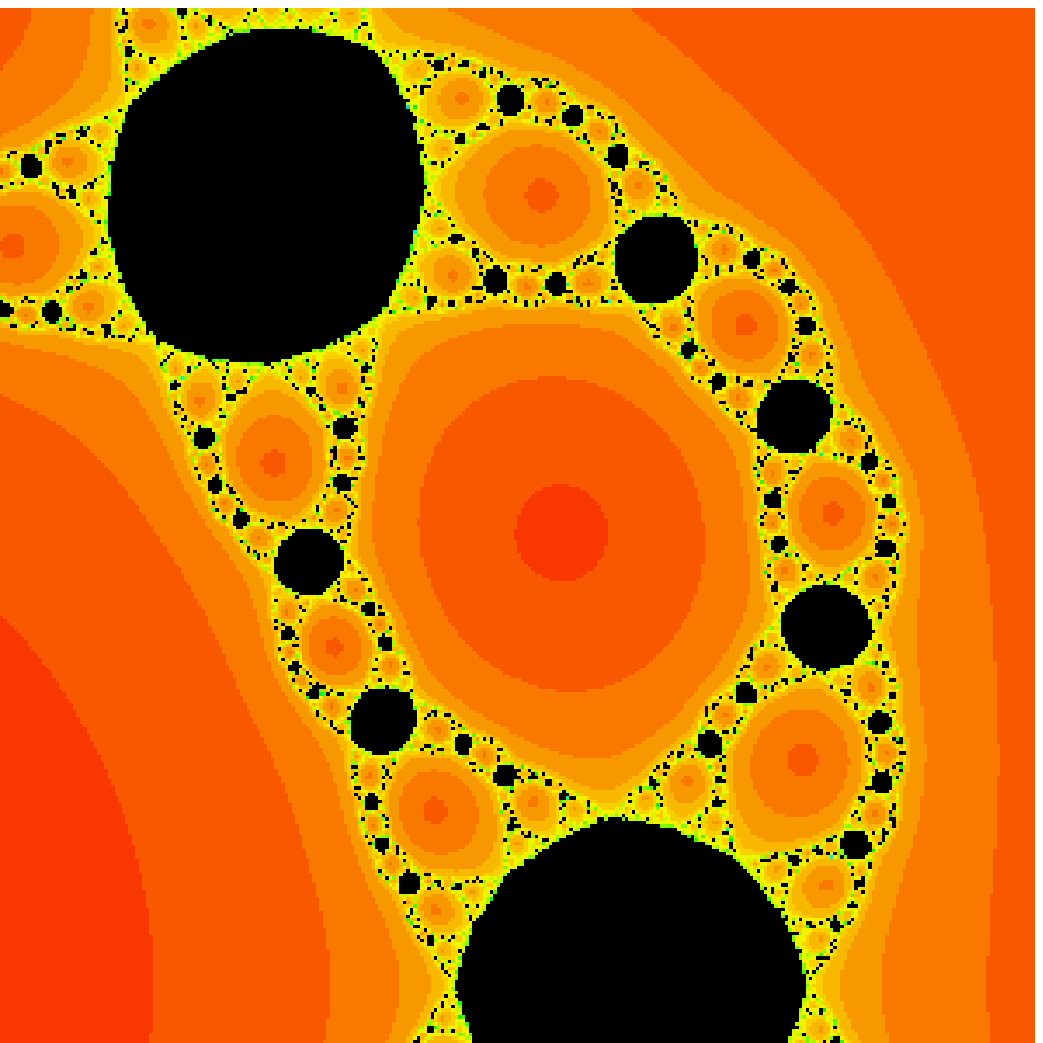}}
\caption{The image on the left is the checkerboard Julia set for 
$F_{0.18}(z)=z^4 + 0.18/z^3$. The image on the right is a 
magnification of one-seventh of the Julia set.}
\label{ckbd}
\end{center}
\end{figure}

Parameters from these cardioids have two distinct types of
Fatou components (see Figure~\ref{ckbd}).  Since $\infty$ is a superattracting fixed point, 
the immediate basin of $\infty$ and its preimages lie in the Fatou set.  These
components are the escaping Fatou components.  The Fatou set also contains
a collection of
components corresponding to 
other finite attracting periodic orbits and their preimages.
These components are the non-escaping
Fatou components.  As we shall show, none of the boundaries of the 
escaping Fatou components
intersect. Likewise, the boundaries of the 
non-escaping Fatou components do not intersect.
However, each such boundary
intersects infinitely many boundaries of the escaping Fatou components, 
and each boundary of an escaping Fatou
component intersects infinitely many non-escaping boundaries. 
Hence, the topology of these Julia sets is very different from the topology
of Sierpi\'nski curve Julia sets.
We use the word ``checkerboard" to describe 
this pattern of Fatou components.

In Figure~\ref{ckbd}, we display the Julia
set for the map $F_{0.18}(z)=z^4 + 0.18/z^3$.
The
red regions are the preimages of the
attracting basin of $\infty$, and the black regions are
the preimages of the basins of the finite attracting cycles.
The boundary of each red region touches infinitely many boundaries of 
the black regions, but it does not touch the boundary of any other red region.
Similarly, the boundary of each black region touches infinitely many boundaries of 
the red regions, but it does not touch the boundary of any other black region.

The external red region in the left-hand image of Figure~\ref{ckbd} is
the immediate basin of attraction of $\infty$. 
We denote this basin by~$\B$. The
central red region that contains the pole at the origin is mapped
to~$\B$. We let $\T$ denote this Fatou component. 
All other red regions are also eventually mapped to~$\B$.

Note that there are $n+d$ ($=7$ in this example)
large  black regions that touch $\B$ and $\T$
at unique points.  These Fatou components are 
(eventually) periodic. We call them the connecting (Fatou)
 components since they are the only Fatou components that extend from
$\B$ to~$\T$. Each of these connecting components 
seems to be separated by another red region
that touches exactly two boundaries of the adjacent connecting components.  
On one side of these red regions, we see $d-1$ ($=2$ in this example)
smaller black
components. On the other side, we see $n-1$ ($=3$ in this example)
smaller black components. Each such black component connects to a pair of red regions.
If we were to magnify this image, we would see that this pattern repeats itself at any scale.

In this paper, we make this construction precise.  In particular,
we give an  algorithm for describing the topological structure of
these Julia sets.  This algorithm also describes the dynamics
on these Julia sets via symbolic dynamics.  We prove:

\bigskip
\noindent{\bf Theorem 1.} {\it Let $\F(z) = z^n + \la/z^d$ with $n
\geq 2$ and $d \geq 1$. Any two Julia sets that correspond to parameters 
in the main cardioids
of the principal Mandelbrot sets in the parameter plane for
these maps are homeomorphic.} 

\medskip
Theorem~1 says that
checkerboard Julia sets are analogous
to Sierpi\`nski curve Julia sets~\cite{W} because 
all checkboard Julia sets are homeomorphic.
 
As in the
Sierpi\`nski case~\cite{DP}, only
certain symmetrically located cardioids give rise to conjugacies on
their respective Julia sets. However, unlike the Sierpi\`nski case, 
we can define 
a
dynamical invariant for checkboard Julia sets. We call it the minimal rotation number, and we prove
that it is a conjugacy invariant for checkerboard Julia sets.

\bigskip
\noindent{\bf Theorem 2.} {\it Two maps drawn from different main
cardioids of principal Mandelbrot sets are topologically conjugate on
their Julia sets if and only if their minimal
rotation numbers are equal.   In particular, two such
maps restricted to their Julia sets
are topologically conjugate only if 
the parameters are 
symmetric either under the rotation $z \mapsto \nu^{j(d+1)} z$ 
or under the map $z \mapsto \nu^{j(d+1)} \overline{z}$,   
where $j \in\bbZ$ and $\nu^{n-1} = 1$. }

\medskip
Theorem~2 leads to an exact count of the number of main
cardioids that have non-conjugate dynamics.

\bigskip
\noindent{\bf Theorem 3.} {\it  Let $g$ be the greatest common divisor
of $n-1$ and $n+d$.  If $g$ is even, then there are exactly $1+g/2$
distinct conjugacy classes among the maps drawn from the main
cardioids of the principal Mandelbrot sets.  If $g$ is odd, then the
number of conjugacy classes is $(g+1)/2$.
}

\bigskip

\section{Preliminaries}

Consider the family of maps
$$
\F(z) = z^n + \frac{\la}{z^d}
$$
where $n \geq 2$, $ d \geq 1$, and $\la \in \bbC$.  
The point at infinity is superattracting of 
order~$n$.  As above, we
denote the immediate basin of $\infty$ by
$\B$.  Also, $0$ is a pole of order~$d$, 
so there is a neighborhood of $0$
that is mapped into $\B$.  If this neighborhood is disjoint
from $\B$, we use the term ``trap door" for the preimage of $\B$ 
that contains $0$. We denote the trap door by  $\T$.

The map $\F(z)$ has $n+d$ ``free" critical points.
They satisfy the equation
$$z^{n+d} = \frac{d\lambda}{n}.$$
Hence, they are equally spaced on the circle of radius 
$$\sqrt[n+d]{\frac{d |\lambda|}{n}}$$ 
centered at the origin.
There are also $n+d$ prepoles. They satisfy the equation
$z^{n+d} = -\lambda$.

The family $\F$ has symmetries in both the dynamical plane 
and the parameter plane. 

\medskip
\noindent\textbf{Symmetry Lemma 1.}\textit{
The map $\F$ is conjugate to $F_{\overline{\lambda}}$ by the conjugacy
$z \mapsto \overline{z}$.
}

\medskip
This first symmetry implies that the parameter plane is symmetric under
complex conjugation.

\medskip
\noindent\textbf{Symmetry Lemma 2.}\textit{
If $\omega$ is a $(n+d)$th root of unity, then
$$\F(\omega z) = \omega^n \F(z).$$}

\vspace{-12pt}
This second symmetry implies that the Julia set of $\F$ is symmetric under the 
map $z\mapsto \omega z$. Similarly, $\B$
and $\T$ possess this $(n+d)\,$-fold symmetry.

Moreover, since the free critical points are arranged symmetrically with respect to
$z \mapsto \omega z$, all of the free critical orbits behave
symmetrically with respect to this rotation. However, it is not necessarily true
that all of these critical orbits behave in the same manner. For example, 
consider the map $F_{0.18}(z)=z^4 + 0.18/z^3$ (see Figure~\ref{ckbd}). 
The orbit of the 
free critical point on the positive real axis 
is asymptotic to a fixed point, but the other six free critical orbits are asymptotic
to a pair of period-three orbits. This symmetry also implies that the basins 
of these attracting orbits are arranged symmetrically as well.

\begin{figure}[ht]
\begin{center}
{\includegraphics[height=2.5in]{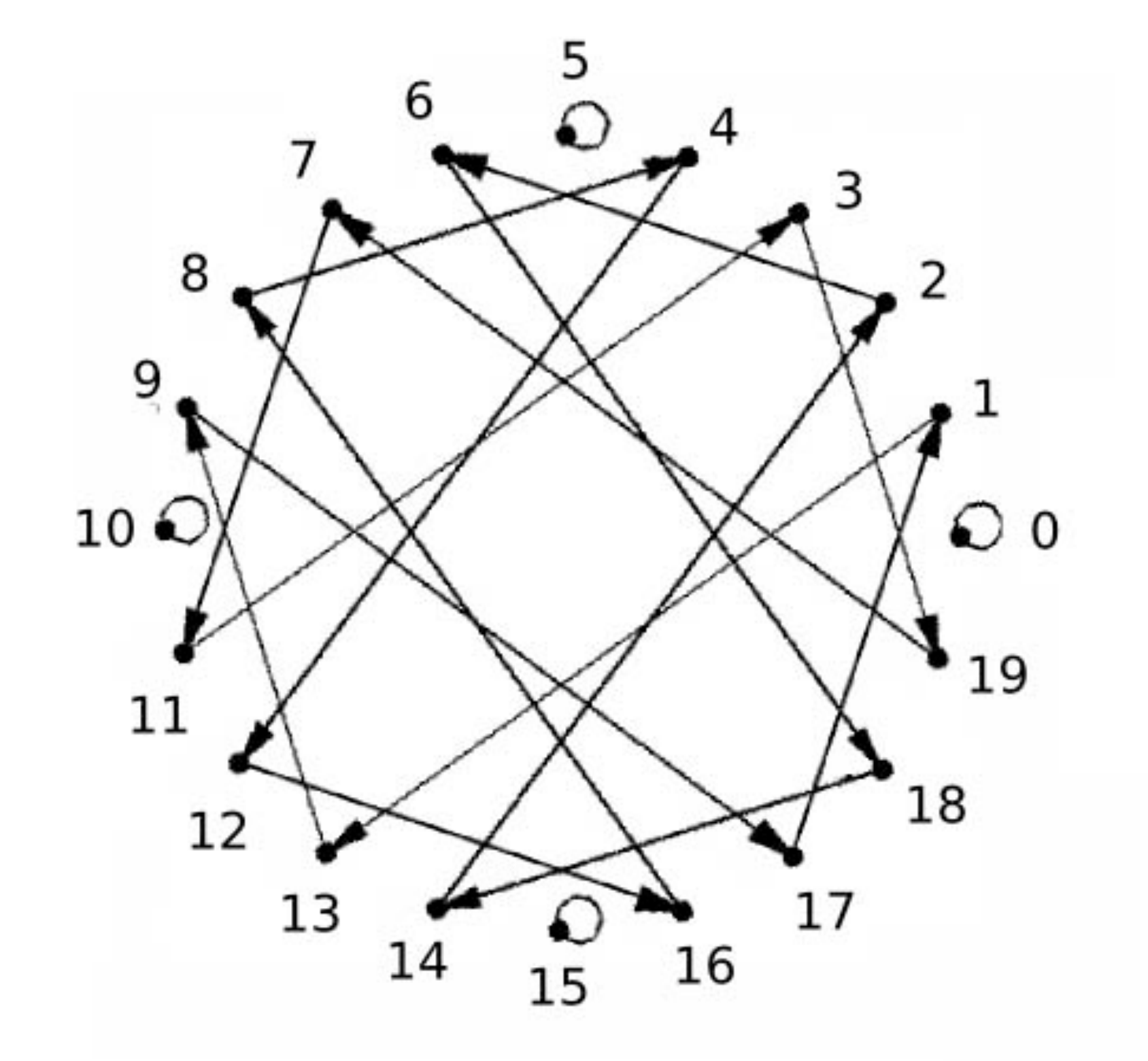}}
\caption{The orbit diagram for the critical points of $F_{\lambda_0}(z)=z^{13} + \lambda_0/z^7$ 
where $\lambda_0$ is the center of the principal Mandelbrot set $\mathcal{M}_0$
that intersects the
positive real axis. The critical point on the positive real axis is a fixed point, and it is labeled with
the number 0. The orbits of the remaining critical points are determined from the orbit
of the fixed point using Symmetry Lemma~2.}
\label{M0n13d7}
\end{center}
\end{figure}

The most important consequence of Symmetry Lemma~2 is the fact that the orbits of all 
of the free critical points can be determined from the orbit of any one of them 
(see Figure~\ref{M0n13d7}). So
the one-dimensional $\la$-plane is a natural parameter plane
for these maps.  

\bigskip
\noindent\textbf{Symmetry Lemma 3.}\textit{
Suppose that $\eta$ is an $(n+d)(n-1)$st root of unity. Let $\nu = \eta^{n+d}$ 
and $\omega=\eta^{n-1}$. Then
$$F^k_{\nu\lambda}(\eta z) = \eta^{n^k}\, F^k_\lambda (z)$$
for $k = 1,\ 2,\ 3, \ldots\,$.}

\medskip\noindent
Note that $\nu$ is an $(n-1)$st root of unity and $\omega$ is
an $(n+d)$th root of unity. Symmetry Lemma~3 is proved by induction on $k$.

\begin{figure}[ht]
\begin{center}
{\includegraphics[height=2.5in]{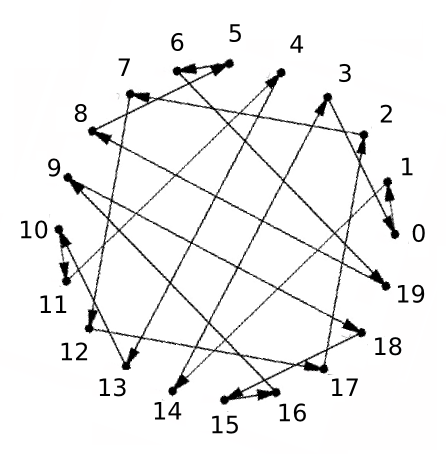}}
\caption{The orbit diagram for the critical points of $F_{\lambda_1}(z)=z^{13} + \lambda_1/z^7$ 
where $\lambda_1 = \nu \lambda_0$ is the center of the principal Mandelbrot set $\mathcal{M}_1$, i.e.,
the image of $\mathcal{M}_0$ under the rotation $z\mapsto \nu z$.}
\label{M1n13d7}
\end{center}
\end{figure}

We can determine the orbit diagram of $F_{\nu\lambda}$ from the orbit diagram of
$\lambda$. In particular, if $c_\lambda$ is a critical point for $\F$, then $\eta c_{\lambda}$ is a critical point
for $F_{\nu\lambda}$. We denote this critical point by $c_{\nu\lambda}$. From Symmetry
Lemma~3, we have
$$
F^k_{\nu\lambda}\left(c_{\nu\lambda}\right)=
F^k_{\nu\lambda}\left(\eta c_{\lambda}\right)=
\eta^{n^k} F^k_{\lambda}\left(c_{\lambda}\right).
$$
Therefore, the orbits of the critical points of $\F$ and $F_{\nu\la}$
behave symmetrically with respect to rotation by some power of $\eta$ (see Figure~\ref{M1n13d7}).
Consequently, 
the parameter plane is symmetric under the rotation
$\la \mapsto \nu \la$.  

As we shall see, the dynamics of $\F$ and
$F_{\nu\la}$ are not necessarily conjugate.  For example, if $\lambda$ lies in 
the main cardioid of the right-hand principal Mandelbrot set in 
the $n=d=3$ case, the map $F_\lambda$ has
a pair of attracting fixed points (see Figure~\ref{pplanes}). In contrast, 
if $\lambda$ is an element of the
the main cardioid of the left-hand principal Mandelbrot set,
the map $F_\lambda$
has
an attracting cycle of period two.

The map
$F_{\nu\lambda}$ is conjugate to $\omega F_\lambda$ by the rotation
$z \mapsto \eta z$ because
$$
F_{\nu\lambda}(\eta z) = 
\left(\eta z\right)^n + \frac{\nu\lambda}{(\eta z)^d} =
\eta^n z^n + \frac{\nu}{\eta^d}\frac{\lambda}{z^d} =
\eta^n\left(z^n + \frac{\lambda}{z^d}\right) =
\eta \omega \left(z^n + \frac{\lambda}{z^d}\right).
$$

More
generally, the map $F_{\nu^j\la}$ is
conjugate to the map $\omega^j\F$
via the rotation $z \mapsto \eta^j z$.

As in the case of the quadratic polynomials $z^2+c$, the orbits of 
the free critical points can tend to $\infty$.  However, unlike the 
quadratic case, there are three distinct ways that these critical orbits 
can escape, and each way leads to a different topological
type of Julia set. 
 
\bigskip
\noindent{\bf Theorem} (The Escape Trichotomy \cite{DLU}.) 
{\it Let $v_\la =
\F(c_\la)$ be a critical value.  
\begin{enumerate}
\item If $v_\la$ lies in $\B$, then $J(\F)$ is a Cantor set.
\item If $\T \neq \B$ and $v_\la$ lies in $\T$, then $J(\F)$ is a Cantor 
set of disjoint, simple, closed curves that surround the origin.
\item In all other cases, $J(\F)$ is a connected set. In particular, 
if $\T\neq \B$ and $\F^j(v_\la)$ lies in $\T$ for some $j \geq 1$, 
then $J(\F)$ is a 
Sierpi\'nski curve. 
\end{enumerate}}
 
We remark that case 2 of this theorem does 
not occur if $n=d=2$
or if $d=1$.
In each parameter plane in Figure~\ref{pplanes},  
the exterior red region consists of
the parameters for which the Julia set is a Cantor set. We call 
this region the 
{\it Cantor set locus}.  The small red region that 
contains the origin is the set of parameters for which 
the Julia set is a Cantor set of circles. We call this region
the {\it McMullen 
domain\/} because McMullen first discovered this type of Julia set 
\cite{McM}.  The complement of these two regions is the 
{\it connectedness locus}. The Julia sets for these parameters are
connected.   The ``holes'' in the connectedness 
locus consist of parameter values for which 
the Julia set is a Sierpi\'nski curve. We call these regions 
{\it Sierpi\'nski holes}. 
 
Note that there are two large Mandelbrot sets
along the real axis 
in the left parameter plane in Figure~\ref{pplanes}.
In the right parameter plane in the same figure,
there are three large Mandelbrot sets symmetrically located 
with respect to the rotation $z\mapsto \nu z$ where
$\nu = e^{2\pi i/3}$.

With the exception of the $n=d=2$ case, 
the parameter plane for the family $\F$ contains
$n-1$ symmetrically located Mandelbrot
sets if $d \neq 1$~\cite{D4}.  We call these sets 
the principal Mandelbrot sets for the family~$\F$. In~\cite{D4},
the existence of these sets was proved 
for the case $n=d>2$. However, the same proof works %
if $d\ne 1$ and $n\neq d$. 
In this paper, we describe the structure of and dynamics on
the Julia sets for parameters that lie in the main
cardioids of these principal Mandelbrot sets.
 
\begin{figure}[ht]
\begin{center}
{\includegraphics[height=2.5in]{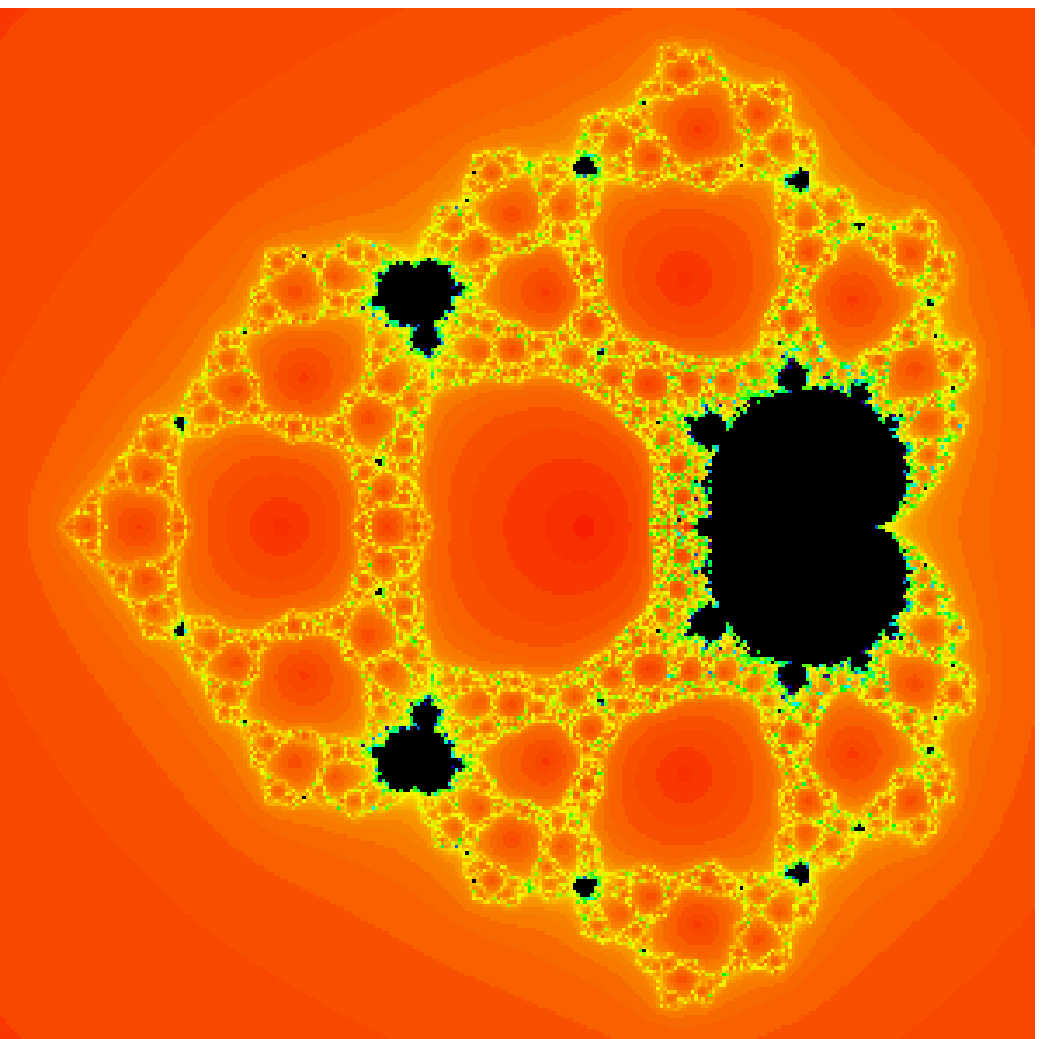}\quad 
\includegraphics[height=2.5in]{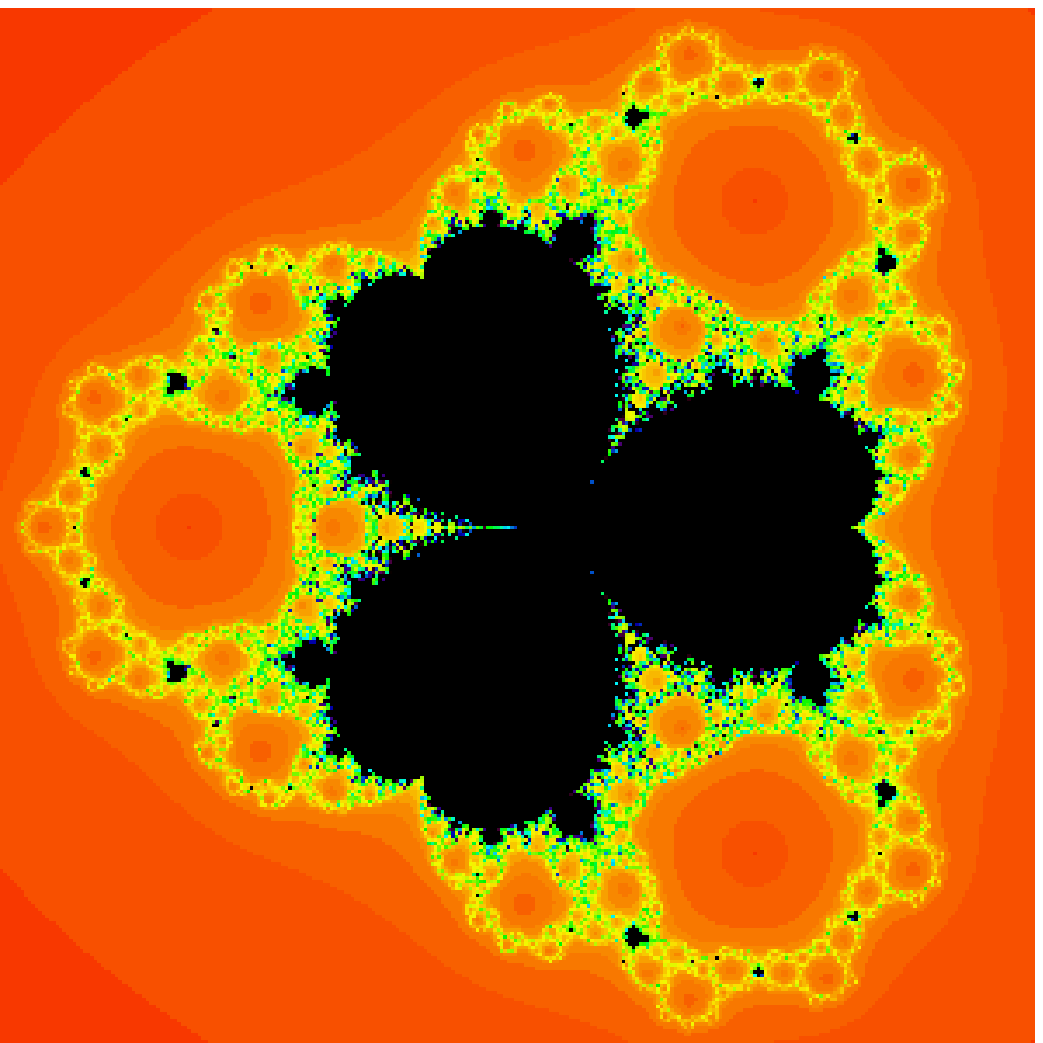}}
\caption{Two parameter planes:  $n=d = 2$ (left)
and $n=4$, $d=1$ (right).}
\label{special}
\end{center}
\end{figure}

If $n=d=2$, there does not exist
a principal Mandelbrot set in the parameter plane.
In this case, the ``tail'' of the Mandelbrot set, i.e., the
parameter corresponding to $c=-2$ in the Mandelbrot set for $z^2 +c$,
extends to the origin, where the map is just $F_0(z)= z^2$. So,
we do not have a
complete Mandelbrot set.  Nonetheless, there is a main cardioid
in which each parameter has 
four connecting Fatou components (see the left-hand parameter
plane in Figure~\ref{special}).

If
$d=1$, 
there are no principal Mandelbrot sets
in the parameter plane.  However, there are $n-1$ distinct
cardioid-shaped regions extending from the Cantor set locus
to the origin (see the right-hand parameter plane in Figure~\ref{special}).  
Parameters drawn from these regions
have $n+1$ connecting Fatou components. 

For the maps that we study, the boundary of the immediate basin
of infinity is a simple, closed curve. 
\smallskip

\noindent{\bf Proposition~1.} {\it  Suppose $\la$ lies in the
connectedness locus in the parameter plane and 
$\F$ is hyperbolic on its
Julia set.  Then the boundary of $\B$
is a simple, closed curve.}

\medskip
\noindent{\bf Proof:}  Since $\F$ is hyperbolic on
its Julia set, $J(\F)$ is locally connected (see~\cite{Mi}).  
In particular, the boundary of $J(\F)$ is locally connected, so
$\partial \B$ is a locally connected set. We need only
show that the set $\overline{\bbC} - \overline{\B}$ is connected.

We argue by contradiction. Suppose that 
$\overline{\bbC} - \overline{\B}$ is disconnected.
Let $W_0$ denote its component 
that contains the trap door.  The second symmetry lemma implies
that 
$W_0$ is symmetric under $z \mapsto \omega z$.  
Also, $\overline{\T}$
is contained in the interior of $W_0$ 
since, if not, there would be a
critical point 
in $\partial \T \cap \partial \B$,
which contradicts the assumption that 
$\F$ is hyperbolic on $J(\F)$.

At least one other
component of $\overline{\bbC} - \overline{\B}$, say $W_1$,  
is mapped over $W_0$. If not, $\partial W_0$ would be
backward invariant, which cannot happen. 
Another application of the second symmetry implies
that $W_j = \omega^j W_1$ is also mapped onto $W_0$
for $j = 1,\ldots,n+d-1$.  We have $n+d$ distinct preimages of
$W_0$.

However,  we claim that there are points in $W_0$ 
that are also mapped
into $W_0$.  To see why, recall that $\partial \T$ is mapped over the
entire boundary of $\B$ by $\F$ and that $\partial \T$ lies in $W_0$.
Thus there is a point $z_0$ in $\partial \T$ that is mapped into
$\partial W_0$.  Then there is a neighborhood of $z_0$ in $W_0$ 
mapped to
a neighborhood of $\F(z_0)$, and hence there are points in this
neighborhood that are mapped inside $W_0$.  
We arrive at a contradiction
since we have found points in $W_0$ 
that have more than $n+d$ preimages.

\hfill{$\Box$}
 
\section{Checkerboard Julia Sets}

In this section we present an algorithm for constructing the Julia sets
for parameters in the main cardioids of the principal Mandelbrot
sets. Theorem~1 follows from this algorithm.

First, consider the case of the principal Mandelbrot set whose
main cardioid intersects the positive real axis, and let $[a,b]$ denote
the interval of intersection.
By considering the graph of $\F\, |\, \bbR$ for $\lambda\in (a,b)$,
we see that each such map
has an attracting fixed point that is
real and positive.  Hence, 
for each $\lambda$ from this main
cardioid, the map $\F$ 
also has an attracting fixed point.  We denote this fixed
point by $p_\la^0$ and its immediate basin of
attraction by $C_\la^0$. Since $\F$ is hyperbolic, 
$\partial C_\la^0$ is a simple, closed curve.  
Furthermore, if $\la \in (a,b)$, 
the graph of $\F\,|\,\bbR$ shows that
$C_\la^0$ extends from $\partial
\B$ to $\partial \T$.  Consequently,
$C_\la^0$ extends from $\partial
\B$ to $\partial \T$ 
for all values of $\lambda$ in this
main cardioid. The intersection  $\partial \B \cap \partial C_\la^0$
consists solely of a repelling fixed point $q_\la^0$ that is real and
positive.
Similarly, the intersection $\partial \T\cap \partial C_\la^0$
is also just one point that is real and positive.
We denote it by 
$u_\la^0$.  Note that $\F(u_\la^0)=q_\la^0$.

From the second symmetry, we obtain
$n+d-1$ other Fatou components that are symmetrically located around the
origin.  We denote
these Fatou components by $C_\la^j$ with $j=1,\ldots,n+d-1$. They are 
ordered in the counterclockwise
direction.
Since each of these components extends from $\T$ to $\B$, we
call them {\it connecting (Fatou) components}.  
Some of these Fatou components are immediate basins of attracting
cycles, and others are eventually periodic components.  The
exact configuration of these components is determined by
Symmetry Lemma~2 with $\omega = \exp(2 \pi i/(n+d))$.  
For example, since $\F(\omega z) = \omega^n
\F(z)$,  we have
$\F(C_\la^1)= C_\la^n$, $\F(C_\la^2)=
C_\la^{2n}$, and so forth.  In particular, if $n=d=3$, both
$C_\la^0$ and $C_\la^3$ are fixed basins, $C_\la^1$ and $C_\la^5$ are
mapped to $C_\la^3$, and $C_\la^2$ and $C_\la^4$ are mapped to
$C_\la^0$.

\begin{figure}[ht]
\begin{center}
{\includegraphics[height=2.5in]{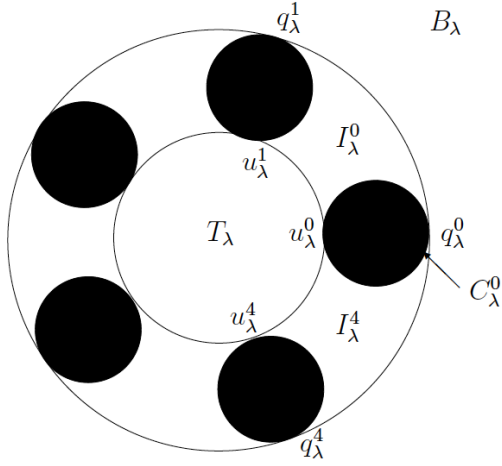}}
\caption{The regions $I_\la^j$ if $n=3$ and $d=2$.}
\label{theIjs}
\end{center}
\end{figure}

Let $p_\la^j = \omega^j p_\la^0$, $q_\la^j= \omega^j q_\la^0$, and
$u_\la^j= \omega^j u_\la^0$.  Then both $q_\la^j$ and $u_\la^j$ lie on $\partial
C_\la^j$. Also, $p_\la^j$ lies in the interior of $C_\la^j$ and is
either periodic or preperiodic (see
Figure~\ref{theIjs}).

For the other $n-2$ main cardioids of the principal Mandelbrot sets,
we have a similar structure due to the $(n-1)\,$-fold symmetry
in the parameter plane.  More precisely, if $\nu$ is an $(n-1)^{\rm st}$
root of unity, the orbits of the critical points of $\F$ and $F_{\nu\la}$ 
behave symmetrically with respect to multiplication by
some (fractional) power of $\nu$, as was shown
immediately following Symmetry Lemma~3.  
Consequently, the configuration of the basins for $F_{\nu\la}$ 
is similar to that of $\F$. 

Recall that $\partial \B$ and $\partial \T$ are simple, closed curves.
Since there are no critical
points in $\partial \B \cap \partial \T$,
these curves do not intersect.
Let $A_\la$ denote the closed annulus bounded by $\partial \B$
and $\partial \T$.
Let $I_\la^j$ denote the closed set in $A_\la$
that is contained in the region located between the open disks 
$C_\la^j$ and $C_\la^{j+1}$. Note that the intersection of $I_\la^j$ and
$I_\la^{j+1}$  is the pair of points $q_\la^{j+1}$
and $u_\la^{j+1}$. Thus there are four points on the boundary of each
$I_\la^j$ that also lie on the boundary of another such set: a pair of
points lies in $I_\la^j \cap I_\la^{j+1}$ and another pair in $I_\la^j
\cap I_\la^{j-1}$. We call the points  $q_\la^j$ the outer junction
points and the points $u_\la^j$ the inner junction points (see
Figure~\ref{theIjs}).   

\bigskip
\noindent{\bf Proposition~2.}
 {\it $\F$ maps each $I_\la^j$ univalently
(except at the junction points) over the region that is the complement
of the three sets $\B$, $\F(C_\la^j)$, and $\F(C_\la^{j+1})$.}

\medskip
\noindent{\bf Proof:} 
 Since $\F$ is conjugate to $z\mapsto z^n$ on
$\partial \B$, the portion of $\partial \B$ that meets $I_\la^j$, 
i.e., the arc in $\partial \B$ connecting $q_\la^{j}$ to $q_\la^{j+1}$,
is mapped to an arc in $\partial \B$ that passes through exactly
$n+1$ outer junction points.  Similarly, the portion of $\partial \T$
that meets $I_\la^j$ is mapped to the complementary arc in $\partial
\B$.  These two arcs meet at a pair of outer junction 
points in $\partial \B$.  Also, the portion of the boundary of
$I_\la^j$ that meets $\partial C_\la^j$ is mapped one-to-one  onto
the boundary of $\F(C_\la^j)$ except
at the junction points. The junction points are both mapped to the same point.
Similarly the other boundary of $I_\la^j$ 
that lies in $\partial C_\la^{j+1}$ is mapped onto
$\partial \F(C_\la^{j+1})$.  Therefore, the boundary of $I_\la^j$ is
mapped to the boundary of the three sets $\B$, $\F(C_\la^j)$, and
$\F(C_\la^{j+1})$.   Since there are no critical points in
$I_\la^j$, the result follows.
\hfill{$\Box$}

\medskip

\begin{figure}[ht]
\begin{center}
{\includegraphics[height=2.5in]{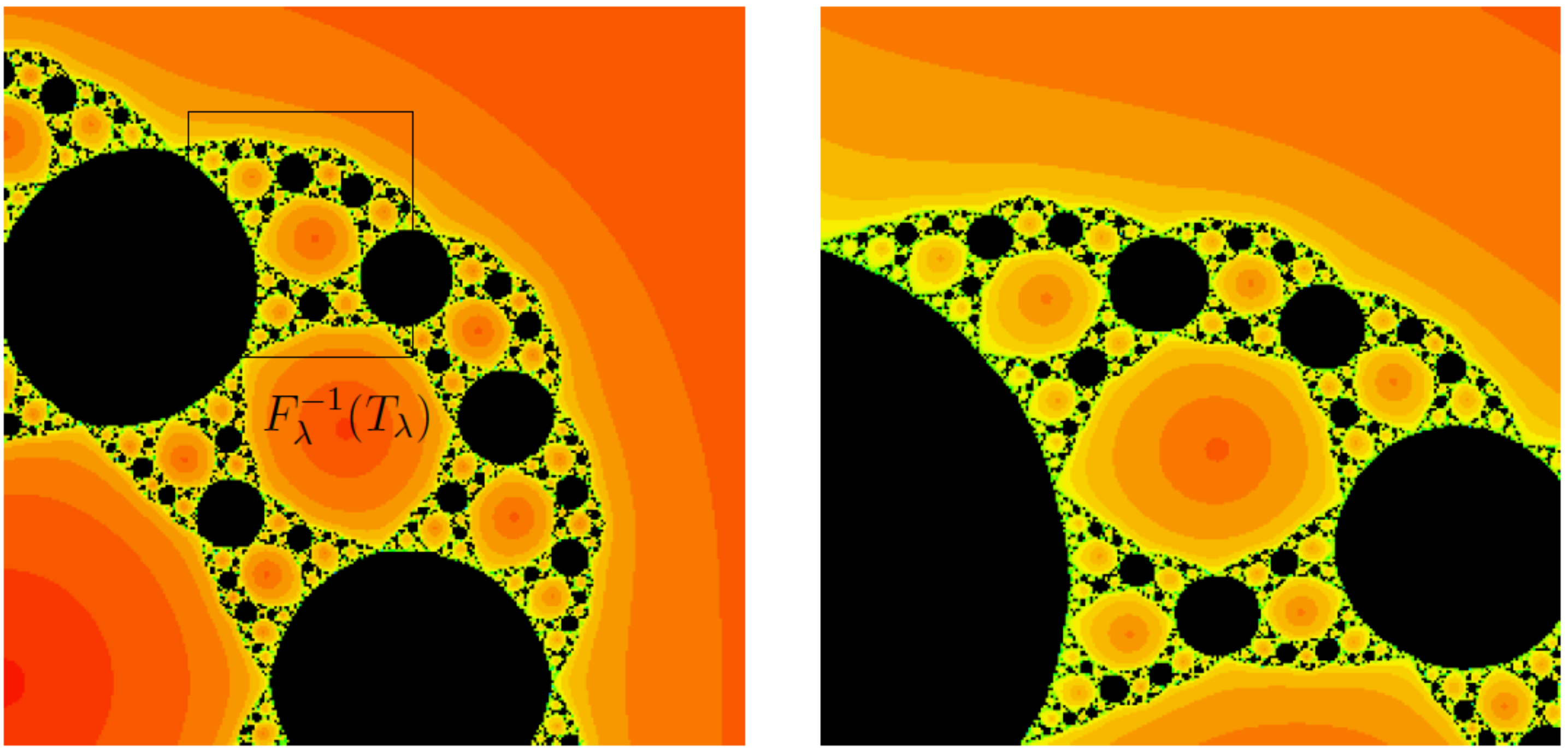}}
\caption{The regions $I_\la^j$ if $n=3$ and $d=2$.}
\label{Ij0-mag}
\end{center}
\end{figure}

We call the two arcs in $I_\la^j$ that lie in the boundaries of
$C_\la^j$ and $C_\la^{j+1}$ the internal
boundary components of $I_\la^j$.
By Proposition~2, there must be a preimage of $\T$
in each $I_\la^j$.  Moreover, the boundary of this preimage must
meet each internal boundary component of $I_\la^j$ in exactly one
point, namely the preimage of the inner junction points lying in the
portions of the boundary of $\F(C_\la^j)$ and $\F(C_\la^{j+1})$ that
lie in $I_\la^j$.
Thus the preimage of $\T$ in each $I_\la^j$ is an open
region whose boundary  meets exactly one point in each of the 
boundaries of the  connecting Fatou components that are adjacent to
$I_\la^j$ (see Figure~\ref{Ij0-mag}).

The preimage of $\T$ separates $I_\la^j$ into two pieces: an
external piece that abuts $\partial \B$ and an internal piece that
abuts $\partial \T$.  The external piece is mapped by $\F$
over the portion of $A_\la$ that stretches from $\F(C_\la^j)$ to
$\F(C_\la^{j+1})$ in the counterclockwise direction.
Since $\F(\omega z) = \omega^n \F(z)$,
this region is mapped over exactly $n$ of
the $I_\la^i$ and $n-1$ of
the $C_\la^i$.  Similarly, the internal piece is mapped over exactly
$d$ of the $I_\la^i$ and 
$d-1$ of the $C_\la^j$. So each of $I_\la^j$ can be further subdivided
as shown in Figure~\ref{Ij0-mag}.  The portion of $I_\la^j$ lying
outside the preimage of $\T$ has $n-1$ preimages of the connecting
components, and the internal portion has $d-1$ such preimages.  Between
each preimage including $C_\la^j$ and $C_\la^{j+1}$, there is a region
that is mapped univalently onto one of the $I_\la^k$'s.  Hence
there is a preimage of each of the sets just constructed in each of these
smaller regions (see Figure~\ref{Ij0-mag}).

Continuing in this fashion, we always find the 
same picture in each region bounded by
$k^{\rm th}$ and earlier preimages of $\T$ and $k^{\rm th}$ and
earlier preimages of the
$C_\la^j$'s. It is a central $(k+1)^{\rm st}$ preimage of $\T$ flanked by
$n-1$ $(k+1)^{\rm st}$  preimages of
the connecting components on one side and $d-1$ other $(k+1)^{\rm st}$ 
preimages on the other side.

However, this construction does not give the entire Julia set of $\F$.  Indeed,
the portion of the Julia set produced thus far contains
only preimages of the boundaries of $\B$ and the $C_\la^j$.  None of these
preimages contain any periodic points; the only periodic points here
lie in $\partial \B$ and $\partial C_\la^j$.  So there must be more to
the Julia set.  

To complete the construction of the Julia set, note that each closed
region $I_\la^j$ is almost mapped univalently over the union of all of the
$I_\la^k$s. The map is univalent 
except at the four junction points. One pair of junction 
points is
mapped to an outer junction point in the image, and the other pair 
is mapped to a
different outer junction point.  We can use symbolic dynamics to
identify each point in the Julia set.  Let $\Sigma$ denote the set of
sequences $(s_0, s_1, s_2, \ldots)$ where each $s_j$ is one of the
integers $0,1,\ldots, n+d-1$.  We identify each point in $J(\F)$ with
a point in $\Sigma$ by assigning to each $z \in J(\F)$ its itinerary
$S(z) = (s_0, s_1, s_2, \ldots) $ where $s_k = j$ if $\F^k(z) \in
I_\la^j$.   However, infinitely many points are assigned to a pair of
sequences. The points $q_\la^j$ and $u_\la^j$ each have a pair of
sequences attached to them since these points reside in two of the
$I_\la^j$'s. For example, 
$S(q_\la^0) = (\overline 0) = 
(\overline{n+d-1})$ and $S(u_\la^0) = (0, \overline{n+d-1}) =
(n+d-1, \overline{0})$.  
Similarly, any point that is eventually
mapped onto a $q_\la^j$ or a $u_\la^j$ also has a pair of itineraries,
e.g., the itineraries $(s_0, \ldots, s_k, 0, \overline{n+d-1})$ and $(s_0,
\ldots, s_k, n+d-1, \overline{0})$ correspond to the same points.  

We
let $\Sigma^\prime$ denote the sequence space with the above
identifications and endow $\Sigma^\prime$ with the quotient topology.
Since each $I_\la^j$ is mapped univalently (except at the junction
points) over the union of the $I_\la^k$ and the Julia set is contained
in this union, standard arguments then show that the Julia set is
homeomorphic to $\Sigma^\prime$.  The subsets $\Sigma_j^\prime$
of $\Sigma^\prime$ consisting of all sequences that start with the
digit~$j$ correspond to points in $I_\la^j \cap J(\F)$, and they 
are homeomorphic to $\Sigma_k^\prime$. It is important to note that the
dynamics on these sets are not the same even though
they are homeomorphic.  
We have described the topological structure of each $I_\la^j \cap
J(\F)$, and this description implies Theorem~1 (see~\cite{Cilingir} for a similar argument).

\section{Dynamical Invariants}
In this section we prove Theorems 2 and 3. Let $\nu=\exp(2 \pi i/(n-1))$. We
show that two maps drawn from the main cardioids of different 
principal Mandelbrot sets are conjugate
on their Julia sets if and only if the cardioids are located
symmetrically under either the maps
$z \mapsto \nu^{j(d+1)}z$ or $z \mapsto
\nu^{j(d+1)}\overline{z}$ for some integer~$j$.

We first observe that it suffices to prove this result for the special
maps whose parameter is the {\it center\/} of these main
cardioids. %
The set of
critical points is invariant under the map, so the critical points are 
either periodic or preperiodic.   
The following proposition follows from the work of Man\'e, Sad, and
Sullivan~\cite{MSS}.

\medskip
\noindent{\bf Proposition~3.} {\it Suppose $\la$ lies at the center of
the main cardioid of one of the principal Mandelbrot sets and that
$\mu$ lies in the same cardioid.  Then $F_\la$ and $F_\mu$ are 
quasiconformally conjugate on their Julia sets.}

\medskip
\noindent{\bf Remark.}  It is not true that $\F$ and $F_\mu$ are
globally conjugate since $\F$ has a superattracting cycle while
the attracting cycle for $F_\mu$ need not be superattracting.

\medskip

By Proposition~3, we need only consider parameters that
lie at the centers of the main cardioids of the principal Mandelbrot
sets.  So for the remainder of this section, we assume that $\la$ and
$\mu$  are centers. %
Then $\mu = \nu^j\la$ 
for some $j \in \bbZ$. The proof of
one direction of Theorem~2 is straightforward:

\medskip
\noindent{\bf Proposition~4.} {\it If $\mu = \nu^{j(d+1)} \la$ or 
$\mu = \nu^{j(d+1)} \overline{\la}$ for some
integer $j$, then $F_\mu$ is conjugate to $\F$.}

\medskip
\noindent{\bf Proof:}  
Let  $\mu = \nu^{j(d+1)} \la$, then 
$$
F_\mu(\nu^j z) = \nu^{jn} z^n + \frac{\la \nu^{j(d+1)}}{\nu^{jd} z^d}
= \nu^j \left( z^n + \frac{\la}{z^d} \right) = \nu^j\F(z).
$$
So $\F$ is  conjugate to $F_\mu$ via the linear map $z \mapsto
\nu^j z$.  

By Symmetry Lemma~1, $\F$ and $F_{\overline{\la}}$ have conjugate
dynamics. So  if $\mu = \nu^{j(d+1)} \overline{\la}$,
then $F_\mu$ is  conjugate to $F_{\overline{\la}}$ and hence
also to $\F$.

\hfill{$\Box$}

From Proposition~4, we know that all centers whose
parameters are of the
form $\nu^k \la$ or $\nu^k\overline{\la}$ where $k =
j(d+1)$ mod$\,(n-1)$ have conjugate
dynamics.   That is, any two main cardioids that are located
symmetrically with respect to either rotation by $\nu^{d+1}$ or
complex conjugation have conjugate dynamics. 
Note that $\nu^{d+1} = \nu^{n-1}\nu^{d+1} =
\nu^{n+d}$, so we can  say that any two cardioids that are
located symmetrically with respect to either rotation by $\nu^{n+d}$
or complex conjugation have conjugate dynamics.

Using basic facts
about the greatest common divisor of two numbers,
we can restate this relationship among the centers with conjugate dynamics
in terms of the greatest common divisor~$g$ of $d+1$ and $n-1$. In fact,
all centers whose parameters are of the form $\nu^k \lambda$ or
$\nu^k \overline{\la}$ where $k$ is an integer multiple of~$g$ have conjugate dynamics.

\medskip
Now we show that these symmetrically located centers are
the only centers with conjugate dynamics.
First we 
define the minimum rotation number for parameters in the main
cardioids of the principal Mandelbrot sets.  For each such parameter
we have $n+d$ connecting components $C_\la^j$ 
with $j$ defined mod $n+d$,
and the $C_\la^j$  are
ordered in the counterclockwise direction as $j$
increases.  Each of these connecting components is mapped two-to-one
onto  another such component since each contains a unique critical
point (see Figure~\ref{theIjs}).   

\begin{figure}[ht]
\begin{center}
{\includegraphics[width=2.5in]{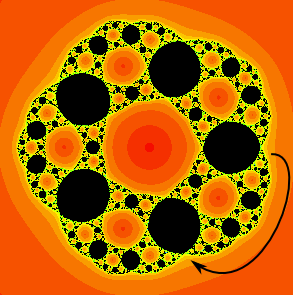}}
\caption{A $C_\la^j$ with $\rho_j = -1$.}
\label{rotationnumber1}
\end{center}
\end{figure}

Suppose $\F(C_\la^j) = C_\la^k$. 
We define the rotation number $\rho_j$ 
of $C_\la^j$ to be the value of
$k-j$ that is closest to $0$ for any $k$ mod $n+d$. Note that it is possible 
for $k$ to be negative (see Figure~\ref{rotationnumber1}).
For example, 
if $\F(C_\la^0) =
C_\la^{n+d-1}$, then the rotation number of $C_\la^0$ would be $-1$
since $C_\la^{n+d-1} = C_\la^{-1}$. %
We say that
$C_\la^j$ is rotated through $k-j$ components if $\rho_j = k$. 
We then
define the {\it minimum rotation number} $\rho(\la)$ 
for $\F$ to be the minimum value of
$|\rho_j|$ over all $j$.  For example, if $\F$ has an attracting fixed
point in some $C_\la^j$, $\rho(\la) = 0$. If there is
no such attracting fixed point, then $\rho(\la) > 0$. 

\medskip
\noindent{\bf Proposition~5.} {\it Let $\la$ and $\mu$ be centers of  main
cardioids of principal Mandelbrot sets.  Then $\F$ is conjugate to
$F_\mu$  if and only if $\rho(\la) = \rho(\mu)$.  Equivalently, $\mu =
\nu^k\la$ or $\mu = \nu^k\overline{\la}$ where $k = j(d+1) \ {\rm mod}\,(n-1)$ for some integer $j$.
The conjugating map
is either %
a rotation about the origin or a rotation followed by complex conjugation.}

\medskip
\noindent{\bf Proof:} 
First suppose
that two such centers $\la$ and $\mu$ have different minimum
rotation numbers.  Then $\F$ and $F_\mu$ cannot be
conjugate on their Julia sets.  To see why, recall that the
connecting components each touch $\partial \B$ at exactly one point.
Now $\partial \B$ must be sent to itself by any conjugacy between $\F$
and $F_\mu$ since this set is the only invariant subset of the Julia
set that touches the boundaries of all of the connecting components.
Thus the ordering of the $C_\la^j$ is either
preserved or reversed by the conjugacy, i.e., either the conjugacy
rotates the connecting components in one direction or the other, or
the conjugacy first applies complex conjugation followed by some
rotation. In either case, if the minimum rotation numbers of $\F$ and
$F_\mu$ are different, then such a conjugacy cannot exist. 

To prove the converse, we consider the 
set of critical points, $c^0_\la, \ldots, c^{n+d-1}_\la$, of $\F$
where $\la$ is the center of a main cardioid of a principal Mandelbrot
set.  The point $c_\la^j$ is the
unique critical point that lies in the connecting component $C_\la^j$.
We note that this set of critical points is invariant under $\F$.

Suppose $\F$ and $F_\mu$ have the same minimum rotation number~$m$. 
By assumption, there is  at least one critical point $c_\la^j$ for
$\F$ for which either $\F(c_\la^j) = c_\la^{j+m} = \omega^m
c_\la^j$ or $\F(c_\la^j) = c_\la^{j-m} = \omega^{-m}
c_\la^j$. (Recall that $\omega^{n+d}=1$.) 
There is also a critical point $c_\mu^i$ for $F_\mu$ for which
either $F_\mu(c_\mu^i) = c_\mu^{i + m} =\omega^m
c_\mu^{i}$ or $F_\mu(c_\mu^i) = c_\mu^{i-m} =
\omega^{-m}c_\mu^{i}$.

We consider the first case for $\la$ and $\mu$, i.e., where the
rotation numbers~$m$ for both critical points are positive.  
Since $\mu = \nu^k\la$ for some $k$,
$F_\mu$ is conjugate to the map $z \mapsto \omega^k \F(z)$
by the rotation $z \mapsto \eta^{k}z$ where
$$\eta=\exp(2\pi i/((n+d)(n-1)))$$ (see the paragraphs 
that follow Symmetry Lemma~3).  So there must be a
critical point for $\omega^k\F$ 
that corresponds to $c_\mu^i$ and
that is also rotated by $\omega^m$ when $\omega^k\F$ is applied to
it. But any critical point of $\omega^k\F$ must also be a critical
point for $\F$. Suppose that $\omega^\ell c_\la^j$ is
the critical point for $\omega^k\F$
that corresponds to $c_\mu^i$.  Then we
have 
$$
\omega^k\F(\omega^\ell c_\la^j)= \omega^m \omega^\ell c_\la^j =
\omega^\ell \omega^m c_\la^j = \omega^\ell \F(c_\la^j).
$$
But
$$
\omega^k\F(\omega^\ell c_\la^j) = \omega^{k+\ell n}\F(c_\la^j).
$$
Therefore we have $\ell = k+\ell n\ {\rm mod}\, (n+d)$.

Consequently,
$\omega^\ell \F(z) = \omega^{k+\ell n}\F(z)$, and using Symmetry
Lemma~2, we obtain
$$\omega^\ell \F(z) =
\omega^k\F(\omega^\ell z)
$$
for all $z \in \bbC$.  
So $\F$ is conjugate to $\omega^k \F$ via the map
$z \mapsto \omega^\ell z$.  Therefore $\F$ is also conjugate to
$F_\mu$ via a linear map of the form $z \mapsto
\eta\omega^\ell z$ and $\mu = \nu^k \la$ where $k =j(d+1)\ {\rm
mod}\, (n-1)$.  

The proof for the case where both rotation numbers~$m$ are negative
is exactly the same.

To prove the second case, the case where the rotations go in opposite directions, we simply conjugate $F_\mu$ to $F_{\overline \mu}$
by complex conjugation and then invoke the first case.  
\hfill{$\Box$}

As a consequence of Propositions 3 and~5,
this result extends to all parameters in the main cardioids of any of
the principal Mandelbrot sets.

\medskip
\noindent{\bf Corollary.} {\it  Let $\la$ and $\mu$ be any parameters
drawn from the main cardioids of any two prinicpal Mandelbrot sets.
Then $\F$ and $F_\mu$ are conjugate on their Julia sets if and only if
$\rho(\la) = \rho(\mu)$.
}

\medskip
Now we can determine exactly which main 
cardioids of the principal 
Mandelbrot sets have conjugate dynamics and the precise number of 
different conjugacy classes. %
 Let $\mathcal{M}_0$ denote the main cardioid of
the principal Mandelbrot set that intersects the positive real axis,
and denote the
remaining main cardioids of principal Mandelbrot sets by $\mathcal{M}_j$ 
where the 
ordering is in the counterclockwise direction.  We write
$\mathcal{M}_j \equiv \mathcal{M}_k$ if the parameters at the centers of
$\mathcal{M}_j$ and $\mathcal{M}_k$ have conjugate dynamics.
Let $g$ be the greatest common divisor of $n-1$
and $d+1$.
As we proved in this section, 
the principal Mandelbrot sets with  dynamics
conjugate to the dynamics of $\mathcal{M}_k$ are those obtained by
successive rotations in the parameter plane by 
$z \mapsto \nu^{jg}z$ or by these rotations followed by complex
conjugation. In particular, we have $\mathcal{M}_0 \equiv
\mathcal{M}_{jg}$ for all integers~$j$. %

\medskip
\noindent {\bf Theorem.} {\it 
If the greatest common divisor $g$ is even, 
there are $1 +g/2$ different conjugacy classes 
among the $\mathcal{M}_j$.  If $g$ is odd, there are
$(g+1)/2$ distinct such conjugacy classes.}

\medskip
\noindent{\bf Proof:}  First suppose that $g=1$.  Then all maps drawn 
from the $\mathcal{M}_j$ have conjugate dynamics, 
so we have $1 = (g+1)/2$
conjugacy classes.  

Now suppose $g>1$.  We claim that $\mathcal{M}_k \not\equiv
\mathcal{M}_0$ for any $k$ with $0 < k < g$.  If not, then maps at the centers of $\mathcal{M}_0$ and some 
$\mathcal{M}_k$ 
would %
be 
conjugate by 
$z \mapsto \overline{z}$
followed possibly by a rotation.  %
But then $\mathcal{M}_k \equiv \mathcal{M}_{-k}$
via $z \mapsto \overline{z}$. 
Also, $\mathcal{M}_g \equiv \mathcal{M}_{-g}$
by $z \mapsto \overline{z}$.  Therefore we have 
$\mathcal{M}_{-g} \equiv \mathcal{M}_{-k}$ by a rotation, %
which would imply that the greatest common divisor 
is smaller than~$g$.  So none of the centers of the $\mathcal{M}_k$ with
$0 < k<g$ have dynamics conjugate to the center of $\mathcal{M}_0$.

If $g$ is even, we
consider $\mathcal{M}_k$ where $0 < k <   g/2$.
We have $\mathcal{M}_k \equiv \mathcal{M}_{-k}$
by complex conjugation. Moreover $\mathcal{M}_{-k} \equiv \mathcal{M}_{g-k}$
since these sets are symmetric under
the rotation $z\mapsto\nu^g z$, so $\mathcal{M}_k \equiv \mathcal{M}_{g-k}$.
On the other hand, we cannot have $\mathcal{M}_k \equiv \mathcal{M}_j$
for any other $j$ with $0 < j<g$ via rotation by $z\mapsto\nu^g z$ or by 
complex conjugation coupled with a rotation, so the principal Mandelbrot
sets with dynamics conjugate to those in $\mathcal{M}_k$ are just
the rotations of $\mathcal{M}_k$ together with their complex conjugates.
The number of such conjugacy classes is $g/2 -1$. 
We have $\mathcal{M}_{g/2} \equiv \mathcal{M}_{-g/2}$ by
the rotation $z\mapsto\nu^{-g} z$ as well as by complex conjugation.  So
$\mathcal{M}_{g/2}$ lies in a conjugacy class that is distinct from the classes of the 
$\mathcal{M}_k$ with $0 \leq k \leq g/2$. The conjugacy class of
$\mathcal{M}_0$ has not yet been counted. Combining all of these classes,
we obtain a total of $1 + g/2$ distinct conjugacy classes.

If $g$ is odd, we count in exactly the same way except that we do not have
a conjugacy class that corresponds to $\mathcal{M}_{g/2}$ in this case.

\hfill${\Box}$

\begin{figure}[ht]
\begin{center}
{
\includegraphics[width=1.8in]{M0_n13d7new.pdf}
\includegraphics[width=1.65in]{M1_n13d7new.png}
\includegraphics[width=1.65in]{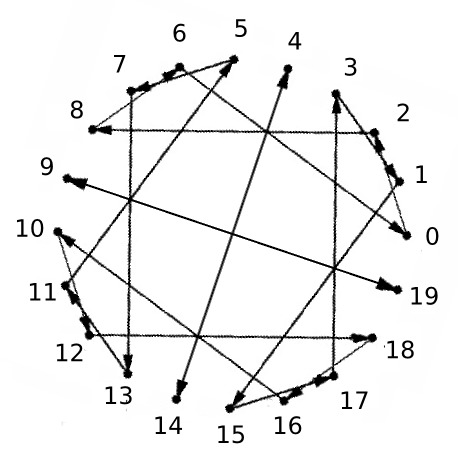}
}
\caption{If $n=13$ and $d=7$, then $g=4$, and consequently, there are three conjugacy classes. 
This figure contains one orbit diagram for each of the three classes.}
\label{threeorbitdiagrams}
\end{center}
\end{figure}

\begin{figure}[ht]
\begin{center}
{
\includegraphics[width=3in]{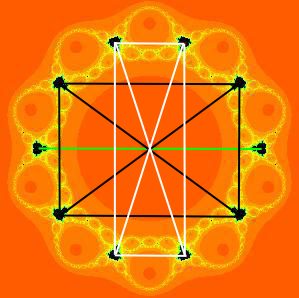}
}
\caption{If $n=11$ and $d=4$, then $g=5$, and consequently, there are three conjugacy classes. 
The parameters with
conjugate dynamics are connected by segments of the same color, e.g., the four Mandelbrot sets
connected by white segments all have conjugate dynamics.}
\label{n11d4conjugaciescolors}
\end{center}
\end{figure}

See Figure~\ref{threeorbitdiagrams} for the three orbit diagrams that arise if $n=13$ and $d=7$. 
In Figure~\ref{n11d4conjugaciescolors}, we consider the case where $n=11$ and $d=4$, and 
group the Mandelbrot sets whose centers have conjugate 
dynamics.

\section{A Group Action}

Since the conjugacies among the $\mathcal {M} _ k$ arise from reflective and rotational symmetries, we can count the number of conjugacy classes by viewing them as orbits of the action of a dihedral group on the set $ \{ \mathcal{M} _ k \} $, viewed as the vertices of a regular $(n-1)$-gon.

Let $a=(n-1)/g$. We claim that the natural group that produces these orbits is $D_{2a}$, the group of symmetries of a regular $a$-gon. Let $s$ be the generator of $D_{2a}$ corresponding to reflection and $r$ be the generator corresponding to rotation. We 
define the action of $D_{2a}$ on $ \{ \mathcal{M} _ k \} $ by 
\begin{eqnarray*}
s \mathcal {M} _ k & = & \mathcal {M} _ {-k \bmod n-1} \\
r \mathcal {M} _ k  & = & \mathcal {M} _ {k+g \bmod n-1}
\end{eqnarray*}

These rules produce 
a well-defined $D_{2a}$ action, and since the actions on $ \{ \mathcal{M} _ k \} $ by $s$ and $r$ are exactly complex conjugation and rotation by $z \mapsto \nu ^ g z$, respectively, the orbits of this action correspond exactly to the conjugacy classes among the $ \mathcal{M} _ k  $. 

By Burnside's Lemma, the number of orbits is 
\[
\frac{1}{|D_{2a}|} \sum_{x \in D_{2a}} | \mathrm {fix} (x)|
 \]
where $\mathrm{fix}(x) = \{ \mathcal{M} _ i \in \{ \mathcal{M} _ k \} : x \mathcal{M} _ i = \mathcal{M} _ i \} $.

The group
$D_{2a}$ has $2a$ elements, and each can be written as $r^j$ or $sr^j$ with $0 \leq j < a$. The identity fixes all $n-1$ elements of $\{ \mathcal{M} _ k \}$, and $r^j$ fixes none for $0<j<a$. Thus the number of orbits is 
\[ 
 \frac{1}{2a} \left( n-1+ \sum_{j=0}^{a-1}| \mathrm {fix} (sr^j)| \right)
\]

An element of the form $sr^j$ rotates each $\mathcal{M} _ k$ by 
$z\mapsto\nu^{jg} z$ and then reflects it about the real axis. 
Equivalently, it reflects the $\mathcal{M} _ k$ through some axis of symmetry of the set viewed as a regular $(n-1)$-gon. Thus, if $n-1$ is odd, every such axis passes through exactly one of the $\mathcal{M} _ k$, and thus $|\mathrm{fix}(sr^j)| = 1$ for all $j$. The formula above then shows the number of conjugacy classes is $ ( n-1+ a)/2a = (g+1)/2$.

If $n-1$ is even, half of the axes of symmetry pass through two of the $  \mathcal{M} _ k $, and half pass through none. Thus $sr^j$ fixes either two or zero of the $\mathcal{M} _ k $. There exists a $j$ such that $sr^j$ fixes none of the $\mathcal{M} _ k$ if and only if there is some $i$ such that $sr^j \mathcal {M} _ i = \mathcal {M} _ {i+1 \bmod n-1}$, i.e., the axis of reflection passes between $\mathcal {M} _ i$ and $ \mathcal {M} _ {i+1 \bmod n-1}$ for some $i$.  For such an $i$, $r^j \mathcal {M} _ i = s^{-1} \mathcal {M} _ {i+1 \bmod n-1} = \mathcal {M} _ {-i-1 \bmod n-1}$ which equals $ \mathcal {M} _ {i+ jg \bmod n-1}$ by the definition of the action of $r$. Thus $-i-1 \equiv i+jg \bmod n-1$,
and hence, $2i + jg +1 \equiv 0 \bmod n-1$. If either $j$ or $g$ is even, 
this equality is impossible since $n-1$ is even, and therefore, $sr^j$ must fix two of the $\mathcal{M} _ k$. If $j$ and $g$ are both odd, however, any $i$ with $ i \equiv (-jg-1)/2 \bmod n-1$ satisfies the congruence, and thus $sr^j$ fixes none of the $\mathcal{M} _ k$. 

Therefore, if $n-1$ and $g$ are even, then $|\mathrm{fix}(sr^j)| = 2$ for all $j$, and the number of conjugacy classes is $( n-1+ 2a)/2a = 1+ g/2$. If $n-1$ is even and $g$ is odd, 
$|\mathrm{fix}(sr^j)|$ equals $2$ if $j$ is even, and $0$ if $j$ is odd. Hence, there are $ ( n-1+ a)/2a = (g+1)/2$ conjugacy classes.

Finally, if $n-1$ is odd, $g$ must be odd, so the possible cases really depend only on the parity of $g$ and not of $n-1$. Hence, the number of conjugacy classes is $(g+1)/2$ if $g$ is odd and $1 + g/2$ if $g$ is even.

\end{document}